\documentclass[12p]{amsart}

\usepackage[english]{babel}
\usepackage{amsmath,amsfonts,amssymb,amscd}
\usepackage[pdftex]{graphicx}

\edef\savecatcodeat{\the\catcode`@}
\catcode`\@=11

\def\tb@ifSpecChars#1#2{#1}
\def\tb@ifNoSpecChars#1#2{#2}

\def\tableau{%
  \bgroup
  \@ifstar{\let\Tif\tb@ifNoSpecChars\tb@tableauB}
          {\let\Tif\tb@ifSpecChars\tb@tableauB}}

\def\tb@tableauB{
  \@ifnextchar[{\tb@tableauC}{\tb@tableauC[]}}

\def\tb@tableauC[#1]{\hbox\bgroup%
    \let\\=\cr
    \def\bl{\global\let\tbcellF\tb@cellNF}%
    \def\tf{\global\let\tbcellF\tb@cellH}
    \dimen2=\ht\strutbox \advance\dimen2 by\dp\strutbox%
    \ifx\baselinestretch\undefined\relax%
    \else%
       \dimen0=100sp \dimen0=\baselinestretch\dimen0%
       \dimen2=100\dimen2 \divide\dimen2 by\dimen0%
    \fi%
    \let\tpos\tb@vcenter
    \tb@initYoung
    \tb@options#1\eoo
    \let\arrow\tb@arrow%
    \dimen0=\Tscale\dimen2%
    \dimen1=\dimen0 \advance\dimen1 by \tb@fframe%
    \lineskip=0pt\baselineskip=0pt
    \def\tb@nothing{}%
    \def\endcellno{$\rss\egroup\bss\egroup}
    \def\endcell{\endcellno\kern-\dimen0}
    \def\begincell{\vbox to\dimen0\bgroup\vss\hbox to\dimen0\bgroup\hss$}%
    \let\overlay\tb@overlay%
    \let\fl\tb@fl%
    \let\lss\hss\let\rss\hss\let\tss\vss\let\bss\vss
    \def\mkcell##1{
        \let\tbcellF\tb@cellD
        \def\tb@cellarg{##1}
        \ifx\tb@cellarg\tb@nothing\let\tb@cellarg\tb@cellE\fi%
        \begincell\tb@cellarg\endcellno
        \tbcellF}
    \let\savecellF\tbcellF
     \Tif{\catcode`,=4\catcode`|=\active}{}\tb@tableauD}%

\let\tb@savetableauD\tableauD
{
    \catcode`|=\active \catcode`*=\active \catcode`~=\active%
    \catcode`@=\active
\gdef\tableauD#1{%
  \Tif{
    \mathcode`|="8000 \mathcode`*="8000%
    \mathcode`~="8000 \mathcode`@="8000%
    \def@{\bullet}%
    \let|\cr
    \let*\tf
    \let~\sk
  }{}%
  \tpos{\tabskip=0pt\halign{&\mkcell{##}\cr#1\crcr}}%
  \global\let\tbcellF\savecellF
  \egroup
  \egroup}
}
\let\tb@tableauD\tableauD
\let\tableauD\tb@savetableauD
\let\tb@savetableauD\undefined

\def\tb@options#1{\ifx#1\eoo\relax\else\tb@option#1\expandafter\tb@options\fi}

\def\tb@option#1{%
  \if#1t\let\tpos\tb@vtop\fi
  \if#1c\let\tpos\tb@vcenter\fi
  \if#1b\let\tpos\vbox\fi
  \if#1F\tb@initFerrers\fi
  \if#1Y\tb@initYoung\fi
  \if#1s\tb@initSmall\fi
  \if#1m\tb@initMedium\fi
  \if#1l\tb@initLarge\fi
  \if#1p\tb@initPartition\fi
  \if#1a\tb@initArrow\fi
}

\def\tb@vcenter#1{\ifmmode\vcenter{#1}\else$\vcenter{#1}$\fi}

\def\tb@vtop#1{\hbox{\raise\ht\strutbox\hbox{\lower\dimen0\vtop{#1}}}}

\def\tb@initPartition{\def\Tscale{.3}}
\def\tb@initSmall{\def\Tscale{1}}
\def\tb@initMedium{\def\Tscale{2}}
\def\tb@initLarge{\def\Tscale{3}}

\def\tb@initArrow{\dimen2=1.25em}

\def\tb@initYoung{%
  \def\tb@cellE{}
  \let\tb@cellD\tb@cellN
  \def\sk{\global\let\tbcellF\tb@cellNF}}
\def\tb@initFerrers{%
  \def\tb@cellE{\bullet}
  \let\tb@cellD\tb@cellNF
  \def\sk{\bullet}}

\tb@initMedium

\def\tb@sframe#1{%
  \vbox to0pt{
    \vss
    \hbox to0pt{%
      \hss
      \vbox to\dimen1{
        \hrule depth #1 height0pt
        \vss
        \hbox to\dimen1{
          \vrule width #1 height\dimen1
          \hss
          \vrule width #1
          }%
        \vss
        \hrule height #1 depth 0in
        }%
      \kern-\tb@hframe
      }%
    \kern-\tb@hframe}}

\def\tb@hframe{.2pt}\def\tb@fframe{.4pt}\def\tb@bframe{2pt}
\def\tb@cellH{\tb@sframe{\tb@bframe}}       
\def\tb@cellNF{}                            
\def\tb@cellN{\tb@sframe{\tb@fframe}}       
\let\tbcellF\tb@cellN                       

\def\tb@rpad{1pt}
\def\tb@lpad{1pt}
\def\tb@tpad{1.8pt}
\def\tb@bpad{1.8pt}

\def\tb@overlay{\endcell\@ifnextchar[{\tb@overlaya}{\begincell}}
\def\tb@overlaya[#1]{\vbox to\dimen0\bgroup%
  \tb@overlayoptions#1\eoo%
  \tss\hbox to\dimen0\bgroup\lss$}
\def\tb@overlayoptions#1{\ifx#1\eoo\relax\else\tb@overlayoption#1\expandafter\tb@overlayoptions\fi}

\def\tb@overlayoption#1{
  \if#1t\def\tss{\vskip\tb@tpad}\let\bss\vss\fi
  \if#1c\let\tss\vss\let\bss\vss\fi
  \if#1b\def\bss{\vskip\tb@bpad}\let\tss\vss\fi
  \if#1l\def\lss{\hskip\tb@lpad}\let\rss\hss\fi
  \if#1m\let\lss\hss\let\rss\hss\fi
  \if#1r\def\rss{\hskip\tb@rpad}\let\lss\hss\fi
}

\def\tb@fl{\endcell\begincell\vrule depth 0pt width \dimen0 height \dimen0 \endcell\begincell}

\@ifundefined{diagram}{}{
\def\tb@arrowpad{.5}

\newoptcommand{\tb@arrow}{\@ne}[2]{%
  \endcell
   \begingroup%
   \let\dg@getnodesize\tb@getnodesize
   \dg@USERSIZE=#1\relax%
   \ifnum\dg@USERSIZE<\@ne \dg@USERSIZE=\@ne \fi%
   \dg@parse{#2}%
   \dg@label{\tb@draw{#1}{#2}}}

\def\tb@getnodesize#1#2#3#4#5{\dimen3=\tb@arrowpad\dimen2 #4=\dimen3 #5=\dimen3\relax}
\def\tb@getnodesize#1#2#3#4#5{\ifnum#2=0\ifnum#3=0\tb@getnodesizetail{#4}{#5}\else\tb@getnodesizehead{#4}{#5}\fi\else\tb@getnodesizehead{#4}{#5}\fi}
\def\tb@getnodesizetail#1#2{\dimen3=.5\dimen2 #1=\dimen3 #2=\dimen3}
\def\tb@getnodesizehead#1#2{\dimen3=.5\dimen2 #1=\dimen3 #2=\dimen3}

\def\tb@draw#1#2#3#4{%
        \dg@X=0\dg@Y=0\dg@XGRID=1\dg@YGRID=1\unitlength=.001\dimen0%
        \dg@LBLOFF=\dgLABELOFFSET \divide\dg@LBLOFF\unitlength%
        \dg@drawcalc
        \begincell
        \let\lams@arrow\tb@lams@arrow
        \begin{picture}(0,0)\begingroup\dg@draw{#1}{#2}{#3}{#4}\end{picture}%
        \endcell
        \endgroup
        \begincell}
}

\def\tb@lams@arrow#1#2{%
 \lams@firstx\z@\lams@firsty\z@
 \lams@lastx#1\relax\lams@lasty#2\relax
 \lams@center\z@
 \N@false\E@false\H@false\V@false
 \ifdim\lams@lastx>\z@\E@true\fi
 \ifdim\lams@lastx=\z@\V@true\fi
 \ifdim\lams@lasty>\z@\N@true\fi
 \ifdim\lams@lasty=\z@\H@true\fi
 \NESW@false
 \ifN@\ifE@\NESW@true\fi\else\ifE@\else\NESW@true\fi\fi
 \ifH@\else\ifV@\else
  \lams@slope
  \ifnum\lams@tani>\lams@tanii
   \lams@ht\ten@\p@\lams@wd\ten@\p@
   \multiply\lams@wd\lams@tanii\divide\lams@wd\lams@tani
  \else
   \lams@wd\ten@\p@\lams@ht\ten@\p@
   \divide\lams@ht\lams@tanii\multiply\lams@ht\lams@tani
  \fi
 \fi\fi
 \ifH@  \lams@harrow
 \else\ifV@ \lams@varrow
 \else \lams@darrow
 \fi\fi
}

\catcode`\@=\savecatcodeat
\let\savecatcodeat\undefined

\newtheorem{theorem}{Theorem}
\newtheorem{proposition}[theorem]{Proposition}

\newtheorem{lemma}[theorem]{Lemma}
\newtheorem{conjecture}[theorem]{Conjecture}
\newtheorem{definition}[theorem]{Definition}
\newtheorem{example}[theorem]{Example}

\title{A non-commutative generalization of $k$-Schur functions}

\author{N. Bergeron \and F. Descouens \and M. Zabrocki }
\address[N. Bergeron]{Dept. of Mathematics and Statistics, York  University, Toronto, Ontario M3J 1P3, CANADA}
\address[F. Descouens]{The Fields Institute, 222 College Street, Toronto, Ontario, M5T 3J1, CANADA}
\address[M. Zabrocki]{Dept. of Mathematics and Statistics, York  University, Toronto, Ontario M3J 1P3, CANADA}

\begin{document}
\maketitle
\def\RR{{\mathbf{R}}}

\begin{abstract}
We introduce non-commutative analogues of $k$-Schur functions of Lapointe-Lascoux and Morse. 
We give an explicit  formulas for the expansions of non-commutive functions with
one and two parameters in terms of these new functions. 
These results are similar to the conjectures existing in the commutative case.
\end{abstract}

\section{Introduction}

In their approach to study Macdonald $q,t$-Kostka polynomials
$K_{\lambda\mu}(q,t)$, L. Lapointe, A. Lascoux and J. Morse
 introduced $k$-Schur functions in \cite{LLM} .  Abstractly, the $k$-Schur functions
 $\{ s_\lambda^{(k)}[X;t] \}$ are generalization of the Schur functions $s_\lambda[X]$
 and are the fundamental basis of
 the subspace of the symmetric functions linearly spanned by the Hall-Littlewood
 symmetric functions 
 $Q'_\lambda[X;t] = \sum_{\mu} K_{\mu\lambda}(0,t) s_\mu[X]$ 
 with all parts of the partition $\lambda$ smaller
 than or equal to $k$.  The motivating property for defining the $k$-Schur functions
 was that the Macdonald polynomials $H_\lambda[X;q,t] = 
 \sum_{\mu} K_{\mu\lambda}(q,t) s_\mu[X]$ 
expand positively in this basis when $\lambda$ has parts bounded by $k$.
 In more recent works \cite{LM1, LM2, LM3, LM4}, Lapointe and Morse have
 studied in particular
 properties of these symmetric functions and 
 conjectured other definitions.  Continuing
 research has established the importance of $k$-Schur functions in other
 areas of mathematics including connections to the geometry of the affine Grassmannian.

In this article we introduce a basis for a subspace of the Hopf algebra
of non-commutative symmetric
functions studied in \cite{NCSF1, NCSF2, NCSF3, NCSF4} 
which we believe is a good analogue of the $k$-Schur functions.  This basis
satisfies many of the same properties of the $k$-Schur functions and,
unlike the commutative counterparts, our analogous versions are
very well behaved so that properties which are difficult to prove or are conjectural
in the commutative case, can be proven for the non-commutative versions.  

Several non-equvialent analogues of the Macdonald and Hall-Littlewood symmetric functions have
been introduced \cite{BZ, Hivert, HLT, Tevlin} to model various
properties of symmetric functions with extra parameters.  In this article we concentrate
on the versions introduced in \cite{BZ} because they have exactly the properties
of the commutative counterparts which we wish to understand better.  In particular,
the non-commutative analogues introduced in \cite{BZ} are known to
have an operator $\blacktriangledown$
which is analogous to the operator $\nabla$ introduced in \cite{BGHT}.  The operator $\nabla$ is defined
so that the Macdonald symmetric functions 
${\tilde H}_\lambda[X;q,t] = t^{n(\lambda)} H_\lambda[X;q,1/t]$ are
eigenfunctions with eigenvalues $t^{n(\lambda)} q^{n(\lambda')}$ 
(here $n(\lambda) = \sum_{i\geq 1} (i-1) \lambda_i$) and similarly $\blacktriangledown$
is defined so that the $q,t$-analogues of non-commutative
symmetric functions are eigenfunctions (see Definition \ref{nabladef}).

The main results that we demonstrate here are difficult to appreciate without
comparing these with the corresponding properties of the $k$-Schur functions.
We list here some of the striking conjectures about the $k$-Schur functions (all of which
except for property \eqref{nablaconj} are from \cite{LLM} or \cite{LM1}) and
preceding each of the statements, we list each of the corresponding theorems
for the non-commutative versions.  Let $\lambda$ be a partition with all parts less
than or equal to $k$.
\begin{enumerate}
\item (Theorem \ref{kSchurBasis})
The $k$-Schur functions are a basis
for the space of symmetric functions 
${\mathcal L} \{ Q'_\mu[X;t] : \mu_1 \leq k \}$.
In particular, the element $Q'_\lambda[X;t]$ expands positively in the
$k$-Schur functions.
\item (Theorem \ref{Branching}) $s_\lambda^{(k)}[X;t]$ expands positively in terms of the functions
$s_\mu^{(k+1)}[X;t]$.
\item (Theorem \ref{MacdonaldPositivity}) 
The Macdonald symmetric function $H_\lambda[X;q,t] = \sum_{\mu} K_{\mu\lambda}(q,t)s_\mu[X]$
expands positively in terms of the $k$-Schur basis
\item (Proposition \ref{omegackSchur} and \ref{omegarkSchur}) 
The fundamental involution $\omega$ applied to $s_{\lambda}^{(k)}[X;t]$ is another
$k$-Schur function with the $t$ parameter inverted.
\item \label{nablaconj} (Theorem \ref{NablakSchurNC})
The operator $\nabla$ acting on $s_\lambda^{(k)}[X;1/t]$ expands positively in the
$k$-Schur basis also with the parameter inverted.
\end{enumerate}

Part of the motivation for considering this filtration of the non-commutative symmetric functions
is that it reflects all of the properties of the $k$-level filtration of the symmetric functions.  
The properties of the $\gamma$-Schur non-commutative functions have motivated
us to study a $k$-level analogue of the $q,t$-Catalan numbers \cite{BDZ} 
and it has also helped make new conjectures on the commutative symmetric functions
through observations
on the non-commutative counterparts (see for instance Conjecture \ref{luckyguess} which
was made after observing the analogous property on the non-commutative versions
of these functions).

\section{Non-commutative symmetric functions}

\subsection{Preliminaries on compositions}
The number of elements on a sequence $\alpha$ is called the length and is denoted by $l(\alpha)$. A sequence of positive integers $\alpha=(\alpha_1, \ldots, \alpha_{l(\alpha)})$ is a composition of size $n$, written $\alpha \models n $, if
\begin{equation*}
\alpha_1 + \ldots + \alpha_{l(\alpha)} = n \ . 
\end{equation*}
A composition $\alpha$ is usually represented by a rim-hook diagram whose rows have lengths $\alpha_1, \ldots, \alpha_{l(\alpha)}$ (read from top to bottom).
\begin{example}\label{ExCompo} The composition $\alpha=(2,4,3,1)$ of size 10 can be represented by the diagram
$$
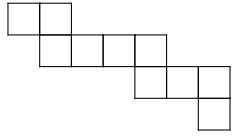
\tableau[sbY]{ & \\ \bl & & & & \\ \bl & \bl & \bl & \bl & & & \\\bl & \bl &\bl & \bl &\bl & \bl &  }
$$
\end{example}
\noindent In the theory of non-commutative symmetric functions, we are  interested in
 two kinds of concatenation
 operations. The first operation is the usual concatenation defined, for two compositions $\alpha$ and $\beta$, by 
\begin{equation}
\alpha \cdot \beta = (\alpha_1, \alpha_2, \ldots, \alpha_{l(\alpha)}, \beta_1, \beta_2, \ldots, \beta_{l(\beta)}) \ .
\end{equation}
The second operation is the attachment defined by
\begin{equation}
\alpha \vert \beta = (\alpha_1, \alpha_2, \ldots,\alpha_{l(\alpha)-1}, \alpha_{l(\alpha)}+\beta_{1},\beta_2, \beta_3, \ldots, \beta_{l(\beta)})\ .
\end{equation}
The descent set $D(\alpha)$ of a composition $\alpha$ is defined as the set
\begin{equation}\label{DescentSet}
D(\alpha) = \{\alpha_1,\ \alpha_1+\alpha_2, \ldots,\alpha_1+\ldots+\alpha_{l(\alpha)-1}\} \ .
\end{equation} 
The descent set $D(\alpha)$ characterizes the composition $\alpha$ and is of size $l(\alpha)-1$. It is easy to see that the compositions of $n$ are in one-to-one correspondence with the subsets of $\{1,2, \ldots, n-1\}$. \\
For any composition $\alpha$, the major index statistic is defined by
\begin{equation}
n(\alpha) = \sum_{i \in D(\alpha)} i = \sum_{i=1}^{l(\alpha)} (i-1)\alpha_{l(\alpha)+1-i}\ .
\end{equation}
For two compositions $\alpha$ and $\beta$, we can refine the previous statistic by defining $c(\alpha,\beta)$ as
\begin{equation}
c(\alpha,\beta)=\sum_{i\in D(\alpha) \cap D(\beta)} i\ .
\end{equation}
There is a natural partial order $\le$ on the set of compositions of $n$, which is called the refinement order. We say that $\alpha$ is finer than $\beta$, written $\alpha \le \beta$, if 
\begin{equation}\label{partialOrder}
D(\beta) \subseteq D(\alpha)\ .
\end{equation}
We can also say that $\alpha \le \beta$ if there exist a sequence of compositions $\gamma^{(1)}, \ldots, \gamma^{(k)}$ such that 
\begin{equation}
\alpha = \gamma^{(1)}\cdot \gamma^{(2)}\cdot \ldots \cdot \gamma^{(k)} \quad \text{and} \quad \beta=\gamma^{(1)}\vert \gamma^{(2)}\vert \ldots \vert \gamma^{(k)} \ .
\end{equation}
There exists three standard involutions on compositions. The first one is the reverse of composition defined by
\begin{equation}
\overleftarrow{\alpha} = (\alpha_{l(\alpha)}, \alpha_{l(\alpha)-1}, \ldots, \alpha_1 ) \ .
\end{equation}
If the descent set of $\alpha$ is $D(\alpha)=\{i_1,\ldots, i_k \}$ then
\begin{equation}
D(\overleftarrow\alpha)=\{ \vert \alpha \vert-i_1, \vert \alpha \vert-i_2, \ldots, \vert \alpha \vert-i_k \} \ .
\end{equation}
The second involution is the complement of a composition. For any composition $\alpha$ of $n$, the complement $\alpha^c$ of $\alpha$ is the composition with descent set the complement of $D(\alpha)$
\begin{equation}
D(\alpha^c)=\{1,2, \ldots, n-1\} \backslash D(\alpha)\ .
\end{equation}
The third one is the analogue of the conjugate of a partition and corresponds to the flipping the composition about the line $y=x$. It is defined by
\begin{equation}
\alpha^\prime\ =\ \overleftarrow{\alpha^c}\ =\ \overleftarrow{\alpha}^c \ .
\end{equation}
\begin{example} The descent set of the composition $\alpha=(2,4,3,1)$ given in Example \ref{ExCompo} is
\begin{equation*}
D(\alpha) = \{2,6,9\}\ .
\end{equation*}
The three previous involutions applied on the composition $\alpha$ give 
\begin{equation*}
\overleftarrow{\alpha} = (1,3,4,2)\ , \ \alpha^c=(1,2,1,1,2,1,2) \ \text{and} \ \alpha^\prime = (2,1,2,1,1,2,1)\ .
\end{equation*}
These compositions correspond respectively to the following diagrams
$$
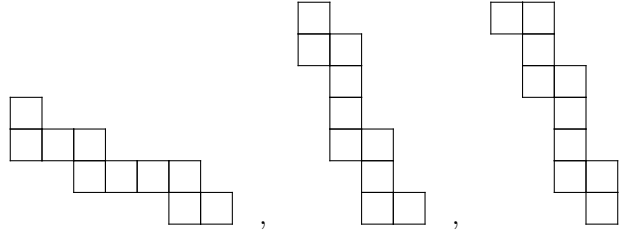
\tableau[sbY]{ \\ & & \\ \bl & \bl & & & & \\ 
               \bl & \bl & \bl & \bl & \bl &  &  }
\quad , \quad
\tableau[sbY]{ \\ & \\ \bl & \\ \bl & \\ \bl & & \\ \bl  & \bl & \\ \bl & \bl & &  }
\quad , \quad 
\tableau[sbY]{ & \\ \bl & \\ \bl & & \\ \bl & \bl & \\ \bl & \bl & \\ \bl & \bl & & 
 \\ \bl & \bl & \bl &} 
$$
\end{example}
\subsection{Non-commutative symmetric functions}
For more details about non-commutative symmetric functions see \cite{NCSF1, NCSF2, NCSF3, NCSF4}. We use the convention of bold font for writting down the non-commutative symmetric functions.
Let $A=\{a_1,a_2,\ldots\}$ be a sequence of non commutative variables and
$X$ the corresponding sequence where variables commute. For any composition
$\alpha$, we define the non-commutative homogeneous functions by
\begin{equation}
\mathbf{h}_\alpha(A) =
\mathbf{h}_{\alpha_1}(A)\ldots\mathbf{h}_{\alpha_{l(\alpha)}}(A),
\end{equation} 
where $\mathbf{h}_{n}(A)$ a non commuting generator of the algebra that
is analogous to the element $h_n(X)$ in the space of symmetric functions.
That is
\begin{equation}
\mathbf{h}_{n}(A) = \sum_{i_1\le i_2\le\cdots\le i_n} a_{i_1}a_{i_2} \cdots
a_{i_n}\,.
\end{equation} 
For any composition $\alpha$, we define the non-commutative homogeneous functions by 
\begin{equation}
\mathbf{h}_\alpha(A) = \mathbf{h}_{\alpha_1}(A)\ldots\mathbf{h}_{\alpha_{l(\alpha)}}(A),
\end{equation}
where $\mathbf{h}_{n}(A)$ a non commuting generator of the algebra that 
is analogous to the element $h_n(X)$ in the space of symmetric functions. 
The product of two non-commutative homogeneous symmetric functions is given by
\begin{equation}
\mathbf{h}_\alpha(A)\mathbf{h}_\beta(A) = \mathbf{h}_{\alpha\cdot\beta}(A) \ .
\end{equation}
The space of non-commutative symmetric functions $\mathbf{Sym}$ over the field $\mathbb{C}(q,t)$ of rational functions in the parameters $q$ and $t$, is defined by 
\begin{equation}
\mathbf{Sym} = \mathbb{C}(q,t) \langle \mathbf{h}_1,\mathbf{h}_2,\ldots \rangle \ .
\end{equation}
The analogs of Schur functions are the ribbon Schur functions defined for any composition $\alpha$ by 
\begin{equation}
\mathbf{R}_\alpha(A) = \sum_{\alpha \le \beta}(-1)^{l(\alpha)-l(\beta)} \mathbf{h}_\beta(A)\ .
\end{equation}
The multiplication rule of two ribbon Schur functions is given by 
\begin{equation}
\mathbf{R}_\alpha(A)\mathbf{R}_\beta(A) = \mathbf{R}_{\alpha\cdot\beta}(A) + \mathbf{R}_{\alpha\vert\beta}(A) \ .
\end{equation}
There are two involutions $\omega^{c}$ and $\overleftarrow \omega$ on the non-commutative symmetric functions which are the analogs of the involution $\omega$ in $Sym$. They are defined for any composition $\alpha$ by
\begin{equation}
\omega^{c} (\mathbf{R}_\alpha(A)) = \mathbf{R}_{\alpha^c}(A) \quad \text{and} \quad \overleftarrow \omega (\mathbf{R}_\alpha(A)) = \mathbf{R}_{\overleftarrow \alpha}(A) \ . 
\end{equation}
We define the commutative evaluation of a non-commutative symmetric function through the surjective map
\begin{equation}
\begin{array}{cccl}
\chi: & \mathbf{Sym}      & \longrightarrow & \mathrm{Sym} \\
      & \mathbf{h}_\alpha(A) & \longmapsto     &  h_{\alpha}(X)
\end{array} \ .
\end{equation}
The image of the ribbon Schur function $\mathbf{R}_\alpha(A)$ by $\chi$ is the commutative skew Schur function indexed by the skew partition corresponding to the ribbon $\alpha$.
\subsection{Deformations of non-commutative symmetric functions}
The modified Hall-Littlewood functions $Q^\prime_\lambda(X;t)$ (resp. modified Macdonald polynomials $\tilde{H}_\lambda(X;q,t)$) are $t$-analogs (resp. $(q,t)$-analogs) of the complete functions $h_\lambda(X)$. In this section, we recall basic statements on the non-commutative analogs of these deformations defined in \cite{BZ}. There exists different non-commutative analogs of Hall-Littlewood functions and Macdonald polynomials which have been considered in \cite{Hivert, HLT} by Hivert, Lascoux and Thibon and more recently in \cite{Tevlin} by Tevlin.

\subsubsection{Non-commutative Hall-Littlewood functions}
In \cite{BZ}, the authors define non-commutative analogs of Hall-Littlewood functions by 
\begin{equation}
\mathbf{H}_\alpha(A;t) = \sum_{\beta \ge \alpha} t^{c(\alpha, \beta^{c})} \mathbf{R}_{\beta}(A) \ . 
\end{equation}
The non-commutative Hall-Littlewood functions $\mathbf{H}_\alpha(A;t)$ satisfy the following specializations
\begin{equation}
\mathbf{H}_\alpha(A;0) = \mathbf{R}_\alpha(A) 
\quad \text{and} \quad 
\mathbf{H}_\alpha(A;1) = \mathbf{h}_\alpha(A) \ .
\end{equation}
\begin{example} The expansion of the non-commutative Hall-Littlewood $\mathbf{H}_{121}(A;t)$ on the ribbon Schur basis is 
\begin{equation*}
\mathbf{H}_{121}(A;t) = \mathbf{R}_{121}(A) + t\ \mathbf{R}_{31}(A) + t^3\ \mathbf{R}_{13}(A) + t^4\ \mathbf{R}_{4}(A) \ .
\end{equation*}
\end{example} 
\noindent For any hook composition $\alpha=(1^a,b)$, the commutative image of $\mathbf{H}_\alpha(A;t)$ by $\chi$ coincides with the commutative modified Hall-Littlewood functions $Q^\prime_{(b,1^a)}(X;t)$
\begin{equation}
\chi(\mathbf{H}_{(1^a,b)}(A;t)) = Q^\prime_{(b,1^a)}(X;t) \ .
\end{equation}
In \cite{BZ}, we find more detailed statements on these non-commutative functions. For example, there are an explicit expansion of the product of two Hall-Littlewood functions in terms of Hall-Littlewood functions, a Pieri formula, some creation operators and a factorization formula at primitive roots of unity. Most of these properties also exist for the analogs considered in \cite{Hivert}.
\subsubsection{Non-commutative Macdonald polynomials}\label{NCMacdo}
In \cite{BZ}, the authors also give a definition for non-commutative analogs of Macdonald polynomials in $\mathbf{Sym}$. These functions are defined by
\begin{equation}\label{MacdoNC}
\mathbf{H}_\alpha(A;q,t) = \sum_{\beta \models\vert \alpha \vert} t^{c(\alpha, \beta^c)}q^{c(\alpha^{\prime}, \overleftarrow\beta)} \mathbf{R}_\beta(A) \ .
\end{equation}
We define non-commutative analogs of the modified Macdonald 
polynomials $\widetilde{\mathbf{H}}_\lambda(X;q,t)$ by
\begin{equation}\label{ModifiedMacdo}
\widetilde{\mathbf{H}}_\alpha(A;q,t) = t^{n(\alpha)} \mathbf{H}_\alpha\left (A;q,\frac{1}{t}\right) = \sum_{\beta \models\vert \alpha \vert} t^{c(\alpha, \beta)}q^{c(\alpha^{\prime}, \overleftarrow\beta)} \mathbf{R}_\beta(A) \ .
\end{equation}
The right hand side of (\ref{ModifiedMacdo}) comes from the following property
\begin{equation}
n(\alpha)-c(\alpha, \beta^c)=\sum_{i \in D(\alpha)}i\quad\quad -\sum_{i\in D(\alpha)\cap(\{1,\ldots, n \}\backslash D(\beta))}i \quad = \sum_{i\in D(\alpha) \cap D(\beta)} i=c(\alpha,\beta) \ .
\end{equation}
\begin{example} The expansion of the non-commutative Macdonald polynomial $\widetilde{\mathbf{H}}_{31}(A;q,t)$ on the ribbon Schur basis is
$$\begin{array}{ccl}
\widetilde{\mathbf{H}}_{31}(A;q,t)& = & \mathbf{R}_{4}(A)\ +\ q^3\ \mathbf{R}_{13}(A)\ +\ q^2\ \mathbf{R}_{22}(A)\ +\ q^5\ \mathbf{R}_{112}(A)\ +\
     t^3\ \mathbf{R}_{31}(A)\ +\\ & &  q^3t^3\ \mathbf{R}_{121}(A)\ +\ q^2t^3\ \mathbf{R}_{211}(A)\ +\ q^5t^3\ \mathbf{R}_{1111}(A) \ . 
\end{array}$$
\end{example}
\noindent These definitions can be expressed in terms of tensor products of $2\times 2$ matrices (see \cite{BZ} for more details). For two matrices $B$ and $C=(c_{ij})_{1\le i,j \le m}$, we use the following convention for the definition of their tensor product
\begin{equation}
B\otimes C = \lbrack c_{ij}B\rbrack_{1 \le i,j\le m}\ .
\end{equation}
Let $n$ be a nonnegative integer. To consider column vectors of elements of degree $n$ in $\mathbf{Sym}$, we need a total order on the set of compositions of $n$. We use the total order corresponding to the rank function $\phi$ defined by 
\begin{equation}
\begin{array}{cccc}
\phi : & \lbrace \alpha, \ \alpha \models n \rbrace & \longrightarrow & \lbrace 0, \ldots, 2^{n-1}-1 \rbrace \\
       & \alpha          & \longmapsto & \sum_{i\in D(\alpha)} 2^{i-1} 
\end{array} \ .
\end{equation}
More precisely, an element of a basis indexed by a composition $\alpha$ appears in the row $\phi(\alpha)$ of the vector.
\begin{example} Using this total order, the compositions of 4 are listed as follows
\begin{equation*}
(4),\ (13),\ (22),\ (112),\ (31),\ (121),\ (211),\ (1111)\ . 
\end{equation*}
\end{example}
\noindent We denote by $\mathbf{R}(A)$ the column vector of ribbon Schur functions $(\mathbf{R}_\alpha(A))_{\alpha \models n}$ and by $\mathbf{H}(A;q,t)$ the column vector of Macdonald polynomials $(\mathbf{H}_\alpha(A;q,t))_{\alpha \models n}$. Directly from Equation (\ref{MacdoNC}), we obtain the following formula  
\begin{equation}
\mathbf{H}(A;q,t) = 
\left \lbrack \begin{array}{cc} 1 & q \\ t^{n-1} & 1 \end{array} \right \rbrack \otimes
\left \lbrack \begin{array}{cc} 1 & q^2 \\ t^{n-2} & 1 \end{array} \right \rbrack  \otimes \ldots \otimes
\left \lbrack \begin{array}{cc} 1 & q^{n-1} \\ t & 1 \end{array} \right \rbrack \mathbf{R}(A)\ .
\end{equation}
The matrix expression for the column vector $\widetilde{\mathbf{H}}(A;q,t)$ of modified non-commutative Macdonald polynomials defined with Equation (\ref{ModifiedMacdo}) is given by
\begin{equation}\label{ModifiedMacdoMatrix}
\widetilde{\mathbf{H}}(A;q,t) = 
\left \lbrack \begin{array}{cc} 1 & q^{n-1} \\ 1 & t \end{array} \right \rbrack \otimes
\left \lbrack \begin{array}{cc} 1 & q^{n-2} \\ 1 & t^2 \end{array} \right \rbrack  \otimes \ldots \otimes
\left \lbrack \begin{array}{cc} 1 & q \\ 1 & t^{n-1} \end{array} \right \rbrack \mathbf{R}(A)\ .
\end{equation}
\begin{example}
The matrix expression of Macdonald polynomials for $n=3$ is given by
\begin{eqnarray*}
\widetilde{\mathbf{H}}(A;q,t)& = &
\left \lbrack \begin{array}{cc} 1 & q^2 \\ 1 & t \end{array} \right \rbrack
\otimes
\left \lbrack \begin{array}{cc} 1 & q \\ 1 & t^2 \end{array} \right \rbrack
\mathbf{R}(A) 
= 
\left \lbrack
\begin{array}{c|cccc}
(3)   & 1 & q^2 & q   & q^3    \\
(12)  & 1 & t   & q   & qt     \\
(21)  & 1 & q^2 & t^2 & q^2t^2 \\
(111) & 1 & t   & t^2 & t^3
\end{array}
\right \rbrack \mathbf{R}(A)\ .
\end{eqnarray*} 
\end{example}
\section{Non-commutative analogs of $k$-Schur functions}
 We define analogues of $k$-Schur functions in the space of non-commutative symmetric functions $\mathbf{Sym}$. To find $k$-Schur functions, Lapointe, Lascoux and Morse originally observed in \cite{LLM} that certain linear combinations of Hall-Littlewood functions were Schur positive and essentially give atoms that make up the Macdonald symmetric functions. \\\\
\noindent Let $n$ be a nonnegative integer and $\gamma$ a composition of $n$. The subspace $\mathbf{Sym}^{(\gamma)}$ of $\mathbf{Sym}$ is defined as the following homogeneous linear span of some non-commutative Hall-Littlewood functions $\mathbf{H}_\alpha(A;t)$
\begin{equation}
\mathbf{Sym}^{(\gamma)} = \mathcal{L} \left \lbrace \mathbf{H}_\alpha(A;t)\ \text{such that} \ \alpha \models \vert \gamma \vert \ \text{and}\ \alpha \le \gamma \right \rbrace \ .
\end{equation}
This space is a natural analogue
 of the subspace of $Sym^{(k)}$ generated by the 
 modified Hall-Littlewood functions $Q^\prime_\lambda(X;t)$
  indexed by partitions $\lambda$ with first part less than $k$. 
\begin{definition}\label{DefGammaSchur}
Let $\alpha$ and $\gamma$ be two compositions of $n$ such that $\alpha \le \gamma$. The $\gamma$-ribbon Schur function $\mathbf{R}_{\alpha}^{(\gamma)}(A;t)$ is defined by 
\begin{equation}
\mathbf{R}_\alpha^{(\gamma)}(A;t) = 
\sum_{\substack{
\beta \ge \alpha\\
D(\alpha)\backslash D(\beta)\subseteq D(\gamma) }} 
t^{c(\alpha, \beta^{c})} \mathbf{R}_\beta(A)\ .
\end{equation}
\end{definition}
\noindent The compositions $\beta$ which appear in the previous sum are those which appear in the interval of the composition poset for the refinement order between $\alpha$ and the composition with descent set $D(\alpha)\backslash D(\gamma)$.
\begin{example} For $n=5$, the expansion of the $(131)$-ribbon Schur function $\mathbf{R}_{1121}^{(131)}(A)$ is 
\begin{equation}
\mathbf{R}^{(131)}_{1121}(A;t) = \mathbf{R}_{1121}(A)+t\ \mathbf{R}_{221}(A)+t^4\ \mathbf{R}_{113}(A)+t^5\ \mathbf{R}_{23}(A) \ .
\end{equation}
\end{example}
\noindent Directly from Definition \ref{DefGammaSchur}, the $\gamma$-ribbon Schur functions reduce for special cases of the level $\gamma$ to some particular non-commutative symmetric functions.
\begin{equation}\label{gamribremark}
\mathbf{R}_{\alpha}^{((\vert \alpha \vert))}(A;t) = \mathbf{R}_\alpha(A) \quad \text{and} \quad \mathbf{R}_\alpha^{(\alpha)}(A;t) = \mathbf{H}_\alpha(A;t)\ . 
\end{equation}
At this moment, it is not clear that the $\gamma$-ribbon Schur functions form a basis of the subspace $\mathbf{Sym}^{(\gamma)}$. The following theorem gives us an explicit expression for the $\gamma$-ribbon Schur functions in terms of Hall-Littlewood functions.
\begin{theorem}\label{kSchurBasis} The set of elements $\left \lbrace \mathbf{R}_\alpha^{(\gamma)}(A;t)\right \rbrace_{\alpha \le \gamma}$ is a basis $\mathbf{Sym}^{(\gamma)}$. The change of bases between Hall-Littlewood functions and $\gamma$-ribbon Schur functions is given by
\begin{equation} \label{HLexpansion}
\mathbf{H}_{\alpha}(A;t) = \sum_{ \alpha \le \beta \le \gamma} t^{c(\alpha,\beta^c)} \mathbf{R}^{(\gamma)}_\beta(A;t)\ .
\end{equation}
The inverse change of basis is given by
\begin{equation}\label{gammaSchurHL}
\mathbf{R}^{(\gamma)}_\alpha(A;t) = \sum_{\alpha \le \beta \le \gamma} (-1)^{l(\alpha)-l(\beta)}t^{c(\alpha, \beta^{c})}\mathbf{H}_\beta(A;t)\ .
\end{equation}
\end{theorem}
\noindent {\bf Proof}: In this theorem, Equation \eqref{gammaSchurHL} is a M\"obius inversion over the
boolean lattice of Equation \eqref{HLexpansion}. The proof of Equation \eqref{HLexpansion} follows from Theorem \ref{Branching} which gives branching rules on $\gamma$-Schur functions and from the limit cases given in Equation \eqref{gamribremark} \hfill $\square$  
\begin{theorem}\label{Branching} Let $\gamma$ and $\widetilde{\gamma}$ be two compositions of $n$ such that $\gamma \le \widetilde{\gamma}$. For any composition $\alpha$ of $n$, the branching rule from the $\gamma$-Schur functions to the $\widetilde{\gamma}$-Schur functions is given by
\begin{equation}
\mathbf{R}_{\alpha}^{(\gamma)}(A;t) = \sum_{\substack{\alpha \le \beta\\ 
D(\alpha) \backslash D(\beta) \subseteq D(\gamma) \backslash D(\widetilde{\gamma})}}t^{c(\alpha,\beta^{c} )}\mathbf{R}_{\beta}^{(\widetilde{\gamma})}(A;t)\ .
\end{equation}
\end{theorem} 
\noindent 
This theorem means that there exists a family of branching rules for a given level $\gamma$. This theorem is also an analogue of the branching rules from the $k$-Schur functions to the $(k+1)$-Schur functions in the commutative case. They are still conjectural and seem to be related to a poset structure of the $k$-shapes. 
\begin{example}  The two branchings for the  $(221)$-Schur function $\mathbf{R}_{(1121)}^{(221)}(A;t)$ to the level $(41)$ and $(23)$ are 
\begin{eqnarray}
\mathbf{R}_{(1121)}^{(221)}(A;t)&=&\mathbf{R}_{(1121)}^{(41)}(A;t)+t^2\ \mathbf{R}_{(221)}^{(41)}(A;t) \\
                                                   &= &\mathbf{R}_{(1121)}^{(23)}(A;t)+t^4\ \mathbf{R}_{(113)}^{(23)}(A;t)\ .
\end{eqnarray}
\end{example}
\noindent In order to prove Theorem \ref{Branching} we need to prove the following two technical lemmas.
\begin{lemma} \label{statequal}
Let $\beta$ be a composition of $n$. For two compositions $\alpha$ and $\delta$ of $n$ such that $\alpha \le \beta$ and $\delta \le \beta$,
$$c(\alpha, \beta^c)+c(\beta, \delta^c) = c(\alpha, \delta^c)~.$$
\end{lemma}  
\noindent
{\bf Proof}:
The quantity $c(\alpha, \beta^c)+c(\beta, \delta^c)$ is the sum over all $i$ in the set $D(\alpha) \cap D(\beta^c)$
and $D(\beta) \cap D(\delta^c)$.  Since $D(\beta) \subseteq D(\alpha)$ and $D(\beta^c)
\subseteq D(\delta^c)$,
\begin{align*}
(D(\alpha) \cap D(\beta^c)) \cup (D(\beta) \cap D(\delta^c)) &=
(D(\alpha) \cap D(\beta^c)) \cup (D(\alpha) \cap D(\beta) \cap D(\delta^c))\\
&= D(\alpha) \cap (D(\beta^c) \cup (D(\beta) \cap D(\delta^c)))\\
&= D(\alpha) \cap ((D(\beta^c) \cap D(\delta^c)) \cup (D(\beta) \cap D(\delta^c)))\\
&= D(\alpha) \cap D(\delta^c) \cap (D(\beta^c)  \cup D(\beta) )\\
&= D(\alpha) \cap D(\delta^c)~.
\end{align*}
Therefore $c(\alpha, \beta^c)+c(\beta, \delta^c) = c(\alpha, \delta^c)$. \hfill $\square$
\begin{lemma} \label{compcor}
\def\tga{{\tilde \gamma}} Let $\alpha$ be a composition of $n$. Let $\gamma$ and $\widetilde{\gamma}$ be two compositions of $n$ such that $\gamma \le \widetilde{\gamma}$. There is a bijection between the two sets
$$A = \{ (\delta, \beta) \text{ such that } \delta \ge \beta \ge \alpha \text{ and } D(\alpha)\backslash D(\beta) \subseteq D(\gamma)\backslash D(\tga) \text{ and } D(\beta)\backslash D(\delta) \subseteq D(\tga) \}$$
and 
$$B = \{ \delta \text{ such that } D(\alpha)\backslash D(\delta) \subseteq D(\gamma) \}~.$$
\end{lemma}
\noindent {\bf Proof}:
\def\tga{{\tilde \gamma}} Let $\alpha$ be a composition of $n$ and assume that $(\delta, \beta)$ is in $A$. Let $i$ be in $D(\alpha) \backslash D(\delta)$, then either $i \in D(\beta)$ or $i \notin D(\beta)$.
If $i \in D(\beta)$, then since $D(\beta) \backslash D(\delta) \subseteq D(\tga)$, 
$$i \in D(\tga) \subseteq D(\gamma)\ .$$  If $i \notin D(\beta)$, then $$i \in D(\alpha) \backslash D(\beta) \subseteq  D(\gamma) \backslash D(\tga) \subseteq D(\gamma)\ .$$
Consequently, $\delta$ is an element of $B$.
\\\\ 
Conversely, consider now the case where $\delta$ is in $B$. Then let
$\beta$ be the composition such that $D(\beta) = D(\delta) \cup D(\tga)$. 
Let $i$ in $D(\alpha) \backslash D(\beta)
= D(\alpha) \backslash (D(\delta) \cup D(\tga))$ such that $i \notin D(\delta)$ and $i \notin D(\tga)$. Hence,
$$i \in D(\alpha) \backslash D(\delta) \subseteq D(\gamma)$$ and therefore $$i \in D(\gamma) \backslash D(\tga)\ .$$
We conclude that $D(\alpha) \backslash D(\beta) \subseteq D(\gamma) \backslash D(\tga)$.\\
Moreover, if $i \in D(\beta) \backslash D(\delta) = (D(\delta) \cup D(\tga)) \backslash D(\delta)$ it must be that $i \in D(\tga)$. Therefore $D(\beta) \backslash D(\delta) \subseteq D(\tga)$.
These two conditions imply that $(\delta, \beta) \in A$.\hfill $\square$
\\\\
{\bf Proof}:(of Theorem~\ref{Branching})
\def\tga{{\tilde \gamma}}
Let $\alpha$ be a composition of $n$ and consider the following expression obtained using Lemma \ref{statequal}
\begin{align}
\sum_{\substack{ \beta \ge \alpha\\
D(\alpha)\backslash D(\beta) \subseteq D(\gamma)\backslash D(\tga)}}  
t^{c(\alpha, \beta^c)}\  \RR_{\beta}^{(\tga)} (A;t) &=
\sum_{\substack{ \beta \ge \alpha\\
D(\alpha)\backslash D(\beta) \subseteq D(\gamma)\backslash D(\tga)}}  
t^{c(\alpha, \beta^c)} \left (
\sum_{\substack{\delta \ge \beta \\ 
D(\beta)\backslash D(\delta) \subseteq D(\tga) } }
t^{c(\beta, \delta^c)}\right ) \ \RR_\delta(A) \nonumber\\
 &=\sum_{\substack{ \beta \ge \alpha\\
D(\alpha)\backslash D(\beta) \subseteq D(\gamma)\backslash D(\tga)}}  
\left ( \sum_{\substack{\delta \ge \beta \\ 
D(\beta)\backslash D(\delta) \subseteq D(\tga) } }
t^{c(\alpha, \delta^c)}\right ) \RR_\delta(A)~.\label{interexpress1}
\end{align}
Lemma \ref{compcor} shows that there is exactly one term in
this sum for every composition in the interval between 
$\alpha$ and the composition with descent set equal to $D(\alpha) \backslash D(\gamma)$
(i.e. compositions $\delta$ such that $D(\alpha) \backslash D(\delta) \subseteq D(\gamma)$).
\\
{}From Definition \ref{DefGammaSchur}, this implies that Equation \eqref{interexpress1} is equal to $\RR^{(\gamma)}_{\beta}(A;t)$. \hfill $\square$
\begin{theorem}\label{MacdonaldPositivity} Let $\alpha$ and $\gamma$ be two compositions of $n$ such that $\alpha \le \gamma$, the non-commutative Macdonald polynomials $\mathbf{H}_\alpha(A;q,t)$ and $\widetilde{\mathbf{H}}_\alpha(A;q,t)$ are $\gamma$-Schur positive. More precisely,
\begin{equation}
\mathbf{H}_\alpha(A;q,t) = \sum_{\beta\le \gamma} t^{c(\alpha,\beta^c)}q^{c(\alpha^\prime, \overleftarrow \beta)} \mathbf{R}_\beta^{(\gamma)}(A) \ ,
\end{equation}
and
\begin{equation}\label{MacdoInkSchur}
\widetilde{\mathbf{H}}_\alpha(A;q,t) = \sum_{\beta\le \gamma} t^{c(\alpha,\beta)}q^{c(\alpha^\prime, \overleftarrow \beta)} \mathbf{R}^{(\gamma)}_\beta\left (A; \frac{1}{t} \right) \ .
\end{equation}
\end{theorem}
\noindent {\bf Proof}:
\def\tga{{\tilde \gamma}}
\newcommand{\invb}[1]{{\overleftarrow{#1}}}
\newcommand{\invc}[1]{{#1^c}}
\newcommand{\inva}[1]{{#1'}}
\def\HH{{\mathbf{H}}}
\begin{align}
\sum_{\beta \le \gamma} t^{c(\alpha,\invc{\beta})} q^{c(\inva{\alpha},\invb{\beta})}
 \RR^{(\gamma)}_{\beta}(A;t)
&=\sum_{\beta \le \gamma} t^{c(\alpha,\invc{\beta})} q^{c(\inva{\alpha},\invb{\beta})} \left (
\sum_{\substack{\delta \ge \beta \nonumber\\ 
D(\beta)\backslash D(\delta) \subseteq D(\gamma) } }
t^{c(\beta, \delta^c)} \right ) \RR_\delta(A) \nonumber\\ 
&=\sum_{\beta \le \gamma} \left ( 
\sum_{\substack{\delta \ge \beta \\ 
D(\beta)\backslash D(\delta) \subseteq D(\gamma) } }
t^{c(\alpha,\invc{\beta})+ c(\beta, \invc{\delta})} q^{c(\inva{\alpha},\invb{\beta}) }\right ) \RR_\delta(A)\ .
\label{macdaddy}
\end{align}
Using Lemma \ref{compcor} with $\alpha = \gamma = (1^n)$ and $\tga = \gamma$, we see that there
is a $1-1$ correspondence between the set
$$\{ (\delta, \beta) \text{ such that } \delta \ge \beta \text{ and } D(\beta) \supseteq 
D(\gamma) \text{ and } D(\beta)\backslash D(\delta) \subseteq D(\gamma) \}$$
and the set of all compositions of $n$. Consequently, each $\beta$ which appears in
this sum is determined from the composition $\delta$ and
its descent set is given by 
$$D(\beta) = D(\gamma) \cup D(\delta)\ .$$
Now we want to show that $c(\alpha,\invc{\beta})+ c(\beta, \invc{\delta}) = c(\alpha, \invc{\delta})$ in the specific case where
$\beta$ does not satisfy the conditions of Lemma \ref{statequal}. 
We must come up with an independent argument.  Note that $D(\gamma) \subseteq D(\alpha)$
and $D(\beta) = D(\gamma) \cup D(\delta)$, hence
\begin{align*}
(D(\alpha) \cap D(\invc{\beta})) \cup (D(\beta) \cap D(\invc{\delta})) &=
(D(\alpha) \cap D(\invc{\gamma}) \cap D(\invc{\delta})) \cup ((D(\gamma) \cup D(\delta)) \cap D(\invc{\delta}))\\
&= (D(\alpha) \cap D(\invc{\gamma}) \cap D(\invc{\delta})) \cup (D(\gamma) \cap D(\invc{\delta}))\\
&= ((D(\alpha) \cap D(\invc{\gamma})) \cup D(\gamma)) \cap D(\invc{\delta})\\
&= D(\alpha) \cap D(\invc{\delta})~.
\end{align*}
Moreover, since $D(\gamma) \subseteq D(\alpha)$, we have  $$D(\inva{\alpha}) \cap 
D(\invb{\gamma}) = \emptyset\ .$$ 
Therefore $$D(\inva{\alpha}) \cap D(\invb{\beta}) =  D(\inva{\alpha}) \cap 
(D(\invb{\gamma}) \cup D(\invb{\delta})) = D(\inva{\alpha}) \cap D(\invb{\delta})\ .$$  
We conclude that $$c(\inva{\alpha},\invb{\beta}) = c(\inva{\alpha},\invb{\delta})\ .$$ Finally, we have that Equation \eqref{macdaddy} is equivalent to 
$$\sum_{\beta \le \gamma} t^{c(\alpha,\invc{\beta})} q^{c(\inva{\alpha},\invb{\beta})}
 \RR^{(\gamma)}_{\beta}(A;t)= \sum_{\delta \models n} t^{c(\alpha, \invc{\delta})} q^{c(\inva{\alpha},\invb{\delta}) } \RR_\delta(A)
= \HH_{\alpha}(A;q,t)\ .$$
The expansion for the modified version $\widetilde{\mathbf{H}}_\alpha(X;t)$ is obtained by using Equation \eqref{ModifiedMacdo} in the previous equation.

\hfill $\square$\\\\
\noindent From now, we need to use a different order on composition than the one use in Section \ref{NCMacdo}.
Given a fix $\gamma\models n$, let $$D(\gamma)=\{i_1,\ldots,i_k\} \text{ and } D(\gamma^c)=\{j_1,\ldots,j_{n-k-1}\}\ ,$$ where $i_1<i_2<\cdots<i_k$ and $j_1<j_2<\ldots<j_
{n-k-1}$. \\ 
Given this, let $\sigma_\gamma$ be the unique permutation of $\{1,2,\ldots ,n-1\}$ defined by 
$$
\begin{array}{lcl}
\sigma_\gamma(i_s)=s      &  \text{ for } 1\le s\le k \ ,\\ 
\sigma_\gamma(j_r)=r+k   &  \text{ for } 1\le r\le n-k-1 \ .
\end{array}
$$ 
We then define for all compositions $\alpha$ of $n$ the rank function
\begin{equation}
\phi_\gamma(\alpha)=\sum_{i\in D(\alpha)} 2^{\sigma_\gamma(i)-1} \ .
\end{equation}
Let denote by $\widetilde{\mathbf{H}}\vert_\gamma(A;q,t)$ the column vector of the modified non-commutative Macdonald polynomials $\widetilde{\mathbf{H}}_\alpha(X;q,t)$ indexed by compositions $\alpha$ such that $\alpha \le \gamma$ ordered using $\phi_\gamma$. The expression of Equation (\ref{MacdoInkSchur}) given in  Theorem \ref{MacdonaldPositivity} can be expressed in terms of $1\times 1$ and $2\times 2$-matrices 
\begin{equation}\label{MacdonaldMatrice}
\widetilde{\mathbf{H}}\vert_\gamma(A;q,t) = \bigotimes_{i \in D(\gamma)} \left \lbrack t^i \right \rbrack \bigotimes_{i \not \in D(\gamma)} 
\left \lbrack 
\begin{array}{cc}
1 & q^{n-i} \\
1  & t^i
\end{array}
\right \rbrack  \mathbf{R}^{(\gamma)}\left (A; \frac{1}{t} \right )\ .
\end{equation}
\begin{proposition}\label{omegackSchur}
Let $\gamma$ be a composition of $n$. For any composition $\alpha \leq \gamma$ 
of $n$, let define $\zeta$ the composition with descent set 
$D(\zeta) = D(\gamma^c) \setminus D(\alpha) \cup D(\gamma)$. 
There exists an analogue of the $k$-conjugation given by
\begin{equation}
\omega^{c}\left ( \mathbf{R}^{(\gamma)}_\alpha(A; t) \right ) 
= t^{n(\gamma)} \mathbf{R}^{(\gamma)}_\zeta\left ( A; \frac{1}{t}\right ) \ .
\end{equation}
\end{proposition}
\begin{proof}
To check the exponent of $t$ is correct, we notice that $D(\zeta) \cap D(\beta) = 
(D(\gamma^c) \setminus D(\alpha) )\cap D(\beta) \cup D(\gamma) \cap D(\beta)
= D(\gamma) \cap D(\beta)$ since $D(\beta) \subseteq D(\alpha)$.  Therefore
$n(\gamma) - c(\gamma, \beta) = c(\gamma, \beta^c)$.  Since we also have
that $D(\gamma) \subseteq D(\alpha)$ and $D(\alpha) \setminus D(\beta) \subseteq D(\gamma)$
then it follows that $n(\gamma) - c(\gamma, \beta) = 
c(\gamma, \beta^c) = c(\alpha, \beta^c)$.  

It is also necessary to verify that $D(\beta) \subseteq D(\alpha)$ and 
$D(\alpha) \setminus D(\beta) \subseteq D(\gamma)$
if and only if $D(\beta^c) \subseteq D(\zeta)$
and $D(\zeta) \setminus D(\beta^c) \subseteq D(\gamma)$.  By noticing
that $D(\beta^c) = D(\gamma) \cap D(\beta^c) \uplus D(\gamma^c) \cap D(\beta^c)$
then $D(\gamma) \cap D(\beta^c) \subseteq D(\gamma)$ and $D(\gamma^c) \cap D(\beta^c)
= D(\gamma^c) \setminus D(\alpha)$ we see that $\beta \geq \zeta$.  This also shows that
$D(\zeta) \setminus D(\beta^c) = D(\gamma) \setminus (D(\gamma) \cap D(\beta^c))$
and is hence a subset of $D(\gamma)$.

$$\omega^{c}\left ( \mathbf{R}^{(\gamma)}_\alpha(A; t) \right ) =
\sum_{\substack{
\beta \ge \alpha\\
D(\alpha)\backslash D(\beta)\subseteq D(\gamma) }} 
t^{c(\alpha, \beta^c)} \mathbf{R}_{\beta^c}(A)
=\sum_{\substack{
\beta \ge \alpha\\
D(\alpha)\backslash D(\beta)\subseteq D(\gamma) }} 
t^{n(\gamma) - c(\gamma, \beta)} \mathbf{R}_{\beta^c}(A)\ .
$$

\end{proof}
\begin{proposition} \label{omegarkSchur}
At $t=1$, we have another analogue of the $k$-conjugation given by
\begin{equation}
\overleftarrow \omega \left ( \mathbf{R}^{(\gamma)}_\alpha(A;1) \right ) = \mathbf{R}^{(\overleftarrow \gamma)}_{\overleftarrow \alpha}(A;1) \ . 
\end{equation}
\end{proposition}
The action of $\overleftarrow \omega$ is a consequence of Proposition \ref{mvomegar} which
follows easily from the definitions and hence we do not provide a proof here.
\section{The non-commutative Nabla operator}
\noindent In the theory of commutative symmetric functions, the operator $\nabla$ is defined as the linear operator which admits the modified Macdonald polynomials $\widetilde{H}_\lambda(X;q,t)$ as eigenvectors for the eigenvalues $t^{n(\lambda)}q^{n(\lambda^\prime)}$. This operator is related to the combinatorics of Dyck paths and to the space of diagonal harmonics \cite{BGHT, GarsiaHaglund, GarsiaHaglund2, GarsiaHaiman, HHLRU, HagCatalan, HagLoehr, LW}. In \cite{BZ}, the authors give a non-commutative analogue 
$\blacktriangledown$ of the operator nabla in the space $\mathbf{Sym}$.
\begin{definition}\label{nabladef}
The non-commutative nabla operator $\blacktriangledown$ is the linear operator defined on the basis of non-commutative modified Macdonald polynomials by  
\begin{equation}
\blacktriangledown\left (\widetilde{\mathbf{H}}_\alpha(A;q,t) \right) = 
t^{n(\alpha)} q^{n(\alpha^\prime)} \widetilde{\mathbf{H}}_\alpha(A;q,t) \ .
\end{equation}
\end{definition}
\noindent This definition can be reformulated in terms of $2\times 2$-matrices as proved in \cite{BZ} by
\begin{equation}\label{NablaMatrice}
\blacktriangledown\left ( \widetilde{\mathbf{H}}(A;q,t) 
\right ) = 
\left \lbrack 
\begin{array}{cc}
q^{n-1} & 0 \\ 0 & t
\end{array}
\right \rbrack
\otimes
\left \lbrack 
\begin{array}{cc}
q^{n-2} & 0 \\ 0 & t^2
\end{array}
\right \rbrack
\otimes \ldots \otimes
\left \lbrack 
\begin{array}{cc}
q & 0 \\ 0 & t^{n-1}
\end{array}
\right \rbrack \widetilde{\mathbf{H}}(A;q,t)\ .
\end{equation}
\begin{proposition}[\cite{BZ}]\label{PropSchurPositivity} For all compositions $\alpha$, the non-commutative functions $\blacktriangledown\left ( \mathbf{R}_\alpha(A)\right )$ is ribbon Schur positive, up to a global sign. More precisely, in terms of matrices, we have
\begin{equation} 
\blacktriangledown \left ( \mathbf{R}(A) \right ) = 
\left \lbrack
\begin{array}{cc}
0 & -q^{n-1} t \\ 1 & (t+q^{n-1})
\end{array}
\right \rbrack \otimes
\left \lbrack
\begin{array}{cc}
0 & -q^{n-2} t^2 \\ 1 & (t^2+q^{n-2})
\end{array}
\right \rbrack\otimes \ldots \otimes
\left \lbrack
\begin{array}{cc}
0 & -q t^{n-1} \\ 1 & (t^{n-1}+q)
\end{array}
\right \rbrack \mathbf{R}(A) \ .
\end{equation}
\end{proposition}
\begin{example} The ribbon Schur expansion of $\blacktriangledown\left (\mathbf{R}_{121} (A) \right )$ is given by
$$
\blacktriangledown \left ( \mathbf{R}_{121}(A) \right ) = -q^2t^2\mathbf{R}_{22}(A) -(q^3t^3+q^2t^5)\mathbf{R}_{211}(A) - (q^5t^2+q^2t^3) \mathbf{R}_{112}(A) - (q^6t^2+q^5t^5+q^3t^3+q^2t^6) \mathbf{R}_{111111}(A) \ .
$$
\end{example}
\noindent In the case of commutative symmetric functions, A. Lascoux gives two conjectures that the commutative symmetric functions $\nabla\left (Q^\prime_\lambda\left (X;\frac{1}{t}\right)\right)$ and $\nabla \left (\omega\left ( Q^\prime_\lambda (X;t)\right)\right )$ are Schur-positive, up to a global sign, for any partition $\lambda$.
\begin{example} The conjectures of Lascoux for the Hall-Littlewood function $Q^\prime_{211}(X;t)$
\begin{equation*}
\nabla\left ( Q_{211}^\prime \left (X;\frac{1}{t}\right ) \right ) = - q t^6  s_{1111}(X) - (q t^5  + q t^4 ) s_{211}(X) - q t^3 s_{31}(X) - q t^4 s_{22}(X) \ ,
\end{equation*}
and 
$$
\begin{array}{cl}
\nabla\left ( \omega( Q_{211}^\prime \left (X;\ t\right ) )\right ) = & (t^9+qt^6)s_{1111} + (t^8+t^7+t^6+qt^5+qt^4) s_{211}(X) + (t^7+t^5+qt^4)s_{22}(X) + \\
                                                                                                            &  (t^6+t^5+t^4+qt^3)s_{31}(X) + t^3s_{4}(X) \ .
\end{array}
$$
\end{example}
\noindent We can generalize the conjectures of Lascoux considering the expansion of the previous functions on the $k$-Schur basis in the parameter $1/t$.
\begin{conjecture}\label{luckyguess}
Let $k$ be a non-negative integer and 
$\lambda$ be a partition of size $n$ such that $\lambda_1 \le k$.  
For any integer $k^\prime$ such that $k^\prime \ge k$, the commutative 
symmetric functions 
$$t^{n(n-1)/2}\nabla\left (Q^\prime_\lambda\left (X;\frac{1}{t}\right)\right) \quad 
\text{ and }
\quad  t^{n(\lambda)+n(n-1)/2}\nabla \left (\omega\left ( Q^\prime_\lambda (X;t)\right)\right )$$ 
are positive, up to a global sign, in the basis of the 
$k^\prime$-Schur functions in the parameter $1/t$.
\end{conjecture}
\begin{example} For $\lambda=(311)$ and $k=3$ and $k=4$, we have the following expansions \\\\
$
\begin{array}{cl}
t^{10} \nabla\left ( Q^\prime_{311}\left ( X; \frac{1}{t}\right )\right ) = & (q^3t^3+q^2)s^{(3)}_{11111}\left ( \frac{1}{t}\right ) + (q^3t^5+q^2t^3+q^2t^2)s^{(3)}_{2111}\left ( \frac{1}{t}\right ) + (q^3t^4+q^2t^3)s^{(3)}_{221}\left ( \frac{1}{t}\right ) + \\ &  q^2t^5s^{(3)}_{311}\left ( \frac{1}{t}\right ) 
\end{array}
$\\\\ and \\\\
$
 \begin{array}{cl}
t^{13} \nabla ( \omega (Q^\prime_{311}(X;t)) = & (q^3t^7+q^2t^6+qt^4+qt^3+1) s^{(3)}_{11111}\left ( \frac{1}{t}\right ) + \\ 
& (q^3t^9+q^2t^8+q^2t^7+qt^7+2qt^6+qt^5+t^4+t^3+t^2) s^{(3}_{2111}\left ( \frac{1}{t}\right ) + \\
& (q^2t^9+q^2t^8+qt^8+2qt^7+qt^6+t^6+t^5+t^4+t^3)s^{(3)}_{221}\left ( \frac{1}{t}\right ) + \\
& (q^2t^9+qt^9+qt^8+t^7+t^6+t^5) s^{(3)}_{311}\left ( \frac{1}{t}\right ) + t^8s^{(3)}_{32}\left ( \frac{1}{t}\right ) \ . 
\end{array}
$
\end{example}
\noindent We prove an analogue of this conjecture in the non-commutative case. We also prove that the functions $\blacktriangledown \mathbf{R}^{(\gamma)}_\alpha(A;t)$ are positive in the $\gamma$-Schur basis. 
\begin{theorem}\label{nablakSchurOnkSchur}
Let $\gamma$ be a composition. For any composition $\alpha$, the functions $\blacktriangledown \mathbf{R}^{(\gamma)}_\alpha\left (A; \frac{1}{t} \right )$ are positive in the $\gamma$-Schur basis in the parameter $1/t$ . More precisely, we have
\begin{equation}
\blacktriangledown \mathbf{R}^{(\gamma)}\left (A; \frac{1}{t} \right ) =  \bigotimes_{i \in D(\gamma)} \left \lbrack t^i \right \rbrack \bigotimes_{i \not \in D(\gamma)}
\left \lbrack 
\begin{array}{cc}
0 & -t^iq^{n-i} \\
1  & t^i+q^{n-i}
\end{array}
\right \rbrack 
\mathbf{R}^{(\gamma)} \left (A;\frac{1}{t}\right) \ .
\end{equation}
It is important to remark here that the order on the vector $\mathbf{R}(A)$ is the one given by $\phi_\gamma$.
\end{theorem}
\noindent {\bf Proof}: Inverting Relation (\ref{MacdonaldMatrice}), we can express
the $\gamma$-ribbon Schur functions in terms of modified Macdonald polynomials as follows
\begin{equation}
\mathbf{R}^{(\gamma)}\left (A; \frac{1}{t} \right ) = \bigotimes_{i \in D(\gamma)} \left \lbrack \frac{1}{t^i} \right \rbrack \bigotimes_{i \not \in D(\gamma)} \frac{1}{t^i-q^{n-i}} 
\left \lbrack 
\begin{array}{cc}
t^i & -q^{n-i} \\
-1  & 1
\end{array}
\right \rbrack \widetilde{\mathbf{H}}\vert_\gamma(A;q,t) \ , 
\end{equation}
where $\widetilde{\mathbf{H}}\vert_\gamma(A;q,t)$ represents the column vector of the $\widetilde{\mathbf{H}}_\alpha(A;q,t)$ for $\alpha \le \gamma$. \\
Applying the linear operator $\blacktriangledown$, we obtain
\begin{equation}\label{7}
\blacktriangledown \mathbf{R}^{(\gamma)}\left (A; \frac{1}{t} \right ) = \bigotimes_{i \in D(\gamma)} \left \lbrack \frac{1}{t^i} \right \rbrack \bigotimes_{i \not \in D(\gamma)} \frac{1}{t^i-q^{n-i}} 
\left \lbrack 
\begin{array}{cc}
t^i & -q^{n-i} \\
-1  & 1
\end{array}
\right \rbrack 
\blacktriangledown \widetilde{\mathbf{H}}\vert_\gamma (A;q,t) \ .
\end{equation}
By definition of the operator $\blacktriangledown$ on $\widetilde{\mathbf{H}}_\alpha(X;q,t)$, we have
\begin{equation}\label{8}
\blacktriangledown \widetilde{\mathbf{H}}\vert_\gamma(A;q,t) = 
\bigotimes_{i \in D(\gamma)} \left \lbrack t^i \right \rbrack
\bigotimes_{i\not \in D(\gamma)} 
\left \lbrack
\begin{array}{cc}
q^{n-i} & 0 \\ 
0  & t^i
\end{array}
\right \rbrack \widetilde{\mathbf{H}}\vert_\gamma(A;q,t)\ .
\end{equation}
Consequently, we have
\begin{equation}
\blacktriangledown \mathbf{R}^{(\gamma)}\left (A; \frac{1}{t} \right ) =  \bigotimes_{i \in D(\gamma)} \left \lbrack 1 \right \rbrack \bigotimes_{i \not \in D(\gamma)} \frac{1}{t^i-q^{n-i}} 
\left \lbrack 
\begin{array}{cc}
t^iq^{n-i} & -t^iq^{n-i} \\
-q^{n-i}  & t^i
\end{array}
\right \rbrack 
\widetilde{\mathbf{H}}\vert_\gamma (A;q,t) \ .
\end{equation}
By using Equation \eqref{MacdonaldMatrice}, we obtain
\begin{equation}
\blacktriangledown \mathbf{R}^{(\gamma)}\left (A; \frac{1}{t} \right ) =  \bigotimes_{i \in D(\gamma)} \left \lbrack t^i \right \rbrack \bigotimes_{i \not \in D(\gamma)}
\left \lbrack 
\begin{array}{cc}
0 & -t^iq^{n-i} \\
1  & t^i+q^{n-i}
\end{array}
\right \rbrack 
\mathbf{R}^{(\gamma)} \left (A;\frac{1}{t}\right) \ .
\end{equation}
\hfill $\square$
\begin{theorem} The image of the non-commutative modified Hall-Littlewood functions by the operator $\blacktriangledown$ is $\gamma$-Schur positive, up to a global sign. More precisely, 
\begin{equation}\label{HLkSchur}
\blacktriangledown\left ( \widetilde{\mathbf{H}}\vert_{\gamma}(A;t)\right )= 
 \bigotimes_{i \in D(\gamma)} \left \lbrack t^{2i} \right \rbrack \bigotimes_{i \not \in D(\gamma)}
\left \lbrack 
\begin{array}{cc}
0 & -t^iq^{n-i} \\
t^i  & t^{2i}
\end{array}
\right \rbrack 
\mathbf{R}^{(\gamma)} \left (A;\frac{1}{t}\right) \ .
\end{equation}
\end{theorem}
\noindent {\bf Proof}: The specialization of Equation \eqref{MacdonaldMatrice} at $q=0$ gives us
\begin{equation}
\widetilde{\mathbf{H}}\vert_\gamma(A;t)= 
 \bigotimes_{i \in D(\gamma)} \left \lbrack t^{i} \right \rbrack \bigotimes_{i \not \in D(\gamma)}
\left \lbrack 
\begin{array}{cc}
1 & 0 \\
1  & t^{i}
\end{array}
\right \rbrack 
\mathbf{R}^{(\gamma)} \left (A;\frac{1}{t}\right) \ .
\end{equation}
Using the result of Theorem \ref{nablakSchurOnkSchur}, we obtain the Equation \eqref{HLkSchur}.
\hfill $\square$
\begin{example} For the non-commutative Hall-Littlewood function $\widetilde{\mathbf{H}}_{21}(A;t)$, we have
\begin{equation}
\blacktriangledown\left(\widetilde{\mathbf{H}}_{121}(A;t) \right )=- q^2  t^6  \mathbf{R}_{22}(A) - q^2  t^9  \mathbf{R}_{211}(A) - q^2t^7 \mathbf{R}_{112}(A) -q^2t^{10} \mathbf{R}_{1111}\ .
\end{equation}
\end{example}
\begin{theorem}\label{NablakSchurNC}
The image of the $\gamma$-Schur functions in the parameter $1/t$ by the operator $\blacktriangledown$ is ribbon Schur positive, up to a global sign. More precisely,  
\begin{equation}
\blacktriangledown \left (\mathbf{R}^{(\gamma)}\left (A;\frac{1}{t}\right)\right ) = \bigotimes_{i\in D(\gamma)} 
\left \lbrack \begin{array}{cc} 1 & t^i \end{array} \right \rbrack \bigotimes_{i \not \in D(\gamma)}
\left \lbrack 
\begin{array}{cc}
0 & -q^{n-i}t^i \\
1 & (t^i+q^{n-i})
\end{array}
\right \rbrack \mathbf{R}(A) \ .
\end{equation}
\end{theorem}
\noindent {\bf Remark}: There are many ways for defining non-commutative analogs of commutative symmetric functions. The fact that $\blacktriangledown$ of these analogs are ribbon-Schur positive, up to a global sign, is an interesting property which is shared with the commutative version as it is conjectured in \cite{BDZ}. On the commutative side, these results permits us to define some generaliations of the $(q,t)$-Catalan numbers.
\section{Multivariate version of the $\gamma$-Schur functions}
In all the previous definitions, it is possible to replace the powers of the parameter $t$ 
by products of the sequence of parameters $t_1,\ldots,t_{n-1}$ and the parameter 
$q$ using the sequence $q_1,\ldots,q_{n-1}$. 
Powers of the parameters $t$ and $q$ are always of the form $c(\alpha, \beta)=\sum_{i\in D(\alpha) \cap D(\beta)} i$, for some compositions $\alpha$ and $\beta$. The multivariate versions permit us
to keep track of  the descents which appear in $c(\alpha,\beta)$.  We reserve the presentation
of these multivariate versions as a side note as these refined results detract from the presentation
of the previous sections.
\begin{definition}
For any composition $\alpha$ of $n$, multivariate non-commutative Hall-Littlewood functions 
is defined by 
\begin{equation} 
\mathbf{H}_\alpha(A;t_1,\ldots,t_{n-1}) = \sum_{\beta \ge \alpha} \left (\prod_{i \in D(\alpha)  \cap D(\beta^{c})}t_i\right ) \mathbf{R}_{\beta}(A) \ . 
\end{equation}
\end{definition}
\noindent These functions are related to the non-commutative Hall-Littlewood functions by the specialization $t_i \rightarrow t^i$
\begin{equation}
\mathbf{H}_\alpha(A;t) = \mathbf{H}_\alpha(A;t,t^2,\ldots,t^{n-1})\ .
\end{equation}
\begin{example}
The expansion of the multivariate non-commutative Hall-Littlewood function $\mathbf{H}_{121}(A;t_1,t_2,t_3)$ on the ribbon Schur basis is 
\begin{equation*}
\mathbf{H}_{121}(A;t_1,t_2,t_3) = \mathbf{R}_{121}(A) + t_1\ \mathbf{R}_{31}(A) + t_3\ \mathbf{R}_{13}(A) + t_1t_3\ \mathbf{R}_{4}(A) \ .
\end{equation*}
\end{example}
\begin{definition}
As for non-commutative Hall-Littlewood functions, we  define a multivariate version for non-commutative Macdonald polynomials by
\begin{equation}\label{MacdoNCMulti}
\mathbf{H}_\alpha(A;q_1,\ldots, q_{n-1},t_1,\ldots, t_{n-1}) = \sum_{\beta \models\vert \alpha \vert} \left ( \prod_{i \in D(\alpha)\cap D(\beta^c)}t_i \prod_{i \in D(\alpha^{\prime})\cap  D(\overleftarrow\beta)} q_i \right ) \mathbf{R}_\beta(A) \ .
\end{equation}
The non-commutative multivariate modified Macdonald polynomials is defined by 
\begin{eqnarray}\label{ModifiedMacdoMulti}
\widetilde{\mathbf{H}}_\alpha(A;q_1,\ldots,q_{n-1},t_1,\ldots, t_{n-1}) = &
\left (\prod_{i \in D(\alpha)} t_i\right ) \mathbf{H}_\alpha\left (A;q_1,\ldots,q_{n-1},1/t_1,\ldots, 1/t_{n-1}\right) \\
 = &
\sum_{\beta \models\vert \alpha \vert} \left ( \prod_{i \in D(\alpha)\cap D(\beta)}t_i\right ) \left ( \prod_{i \in D(\alpha^{\prime}) \cap D(\overleftarrow\beta)} q_i \right ) \mathbf{R}_\beta(A) \ .
\end{eqnarray}
\end{definition}
\noindent The non-commutative Macdonald polynomials $\mathbf{H}_\alpha(X;q,t)$ (resp. $\mathbf{H}_\alpha(X, q_1,\ldots,q_n,t_1,\ldots,t_n)$) and their modified versions $\widetilde{\mathbf{H}}_\alpha(X;q,t)$ (resp. $\widetilde{\mathbf{H}}_\alpha(X; q_1,\ldots,q_n,t_1,\ldots,t_n)$) coincide under the specialization $t_i \rightarrow t^i$ and $q_i \rightarrow q^i$.  
\begin{example}
The expansion of the multivariate non-commutative modified Macdonald polynomial \\ $\widetilde{\mathbf{H}}_{31}(A;q_1,q_2,q_3,t_1,t_2,t_3)$ on the ribbon Schur basis is
$$\begin{array}{ccl}
\widetilde{\mathbf{H}}_{31}(A;q_1,q_2,q_3,t_1,t_2,t_3)& = & \mathbf{R}_{4}(A)\ +\ q_3\ \mathbf{R}_{13}(A)\ +\ q_2\ \mathbf{R}_{22}(A)\ +\ q_2q_3\ \mathbf{R}_{112}(A)\ +\
     t_3\ \mathbf{R}_{31}(A)\ +\\ & &  q_3t_3\ \mathbf{R}_{121}(A)\ +\ q_2t_3\ \mathbf{R}_{211}(A)\ +\ q_2q_3t_3\ \mathbf{R}_{1111}(A) \ . 
\end{array}$$
\end{example}
\noindent
For the multivariate non-commutative Macdonald polynomials and their modified version, matricial expressions are given by 
\begin{equation}
\mathbf{H}(A;q,t) = 
\left \lbrack \begin{array}{cc} 1 & q_1 \\ t_{n-1} & 1 \end{array} \right \rbrack \otimes
\left \lbrack \begin{array}{cc} 1 & q_2 \\ t_{n-2} & 1 \end{array} \right \rbrack  \otimes \ldots \otimes
\left \lbrack \begin{array}{cc} 1 & q_{n-1} \\ t_1 & 1 \end{array} \right \rbrack \mathbf{R}(A)\ ,
\end{equation}
\begin{equation}
\widetilde{\mathbf{H}}(A;q,t) = 
\left \lbrack \begin{array}{cc} 1 & q_{n-1} \\ 1 & t \end{array} \right \rbrack \otimes
\left \lbrack \begin{array}{cc} 1 & q_{n-2} \\ 1 & t_2 \end{array} \right \rbrack  \otimes \ldots \otimes
\left \lbrack \begin{array}{cc} 1 & q_1 \\ 1 & t_{n-1} \end{array} \right \rbrack \mathbf{R}(A)\ .
\end{equation}
\begin{definition}
Let $\alpha$ and $\gamma$ be two compositions of $n$ such that $\alpha \le \gamma$.  A multivariate version of non-commutative $\gamma$-Schur functions is defined by
\begin{equation}
\mathbf{R}_\alpha^{(\gamma)}(A;t_1,\ldots,t_{n-1}) = 
\sum_{\substack{ \beta \ge \alpha\\
D(\alpha)\backslash D(\beta)\subseteq D(\gamma)}} 
\left (\prod_{i \in D(\alpha) \cap D(\beta^{c})}t_i\right ) \mathbf{R}_\beta(A)\ .
\end{equation}
\end{definition}
All the results stated in the previous sections can be generalized to the
multivariate versions for the most part simply by changing $t^i \rightarrow t_i$.
For practical notational purposes it was convenient to state the results
using only the two parameters $q$ and $t$.
\begin{proposition}\label{mvomegar}
 The action of the analogue
  of the $k$-conjugation on the multivariate $\gamma$-Schur functions is given by
\begin{equation}
\overleftarrow \omega \left ( \mathbf{R}^{(\gamma)}_\alpha(A;t_1,\ldots,t_{n-1}) \right ) = \mathbf{R}^{(\overleftarrow \gamma)}_{\overleftarrow \alpha}(A;t_{n-1},\ldots,t_1) \ . 
\end{equation}
\end{proposition}
\noindent \textbf{Remark}: The non-commutative analogs of Hall-Littlewood functions and Macdonald polynomials defined in \cite{Hivert} admit also a multivariate version defined in \cite{HLT}.
\\\\
\noindent \textbf{Acknowlegment}
All the computations related to this work have been done using the combinatorics package MuPAD-Combinat for the algebra computer system MuPAD (see \cite{HT} for an introduction to this package and especially \cite{D} for computations and implementations of symmetric functions)

\section{Appendix}
\subsection{Tables of $\gamma$-Schur functions on ribbon Schur functions for weight 4}
These tables are calculated from Definition \ref{DefGammaSchur} .  The columns of the table
index the corresponding elements of the $\gamma$-ribbon Schur basis and the rows index
the subscripts of the ribbon Schur basis.
$$
\begin{array}{|c|cccc|} \hline
          & (31) & (121) & (211) & (1111) \\ \hline
(4)      &  t^3 &   \cdot       &  \cdot     & \cdot  \\
(13)    &  \cdot      & t^3    &  \cdot     & \cdot  \\
(22)    &   \cdot    &  \cdot        &  t^3&  \cdot  \\
(112)  &  \cdot      & \cdot         &   \cdot    &   t^3\\
(31)    &   1   &    \cdot      &   \cdot    &   \cdot   \\
(121)  & \cdot       &  1      &    \cdot   & \cdot    \\
(211)  & \cdot       &  \cdot        & 1    & \cdot  \\
(1111) & \cdot      &    \cdot      &  \cdot     &  1 \\ \hline
\multicolumn{5}{c}{} \\
\multicolumn{5}{c}{(31)\text{-Schur functions}} \\
\end{array} \quad \quad \quad
\begin{array}{|c|cccc|}
\hline 
          & (22) & (112) & (211) & (1111) \\ \hline
(4)      &  t^2 &   \cdot         &  \cdot     & \cdot  \\
(13)    &  \cdot      & t^2      &  \cdot     & \cdot  \\
(22)    &   1    &  \cdot         &  \cdot     &  \cdot  \\
(112)  &  \cdot      & 1           &    \cdot    &   \cdot \\
(31)    &   \cdot  &    \cdot      &   t^2       &   \cdot   \\
(121)  & \cdot       &  \cdot      &    \cdot   & t^2    \\
(211)  & \cdot       &  \cdot        & 1          & \cdot  \\
(1111) & \cdot      &    \cdot      &  \cdot     &  1 \\ \hline
\multicolumn{5}{c}{} \\
\multicolumn{5}{c}{(22)\text{-Schur functions}} \\
\end{array} $$
$$
\begin{array}{|c|cccc|}
\hline
          & (13) & (112) & (121) & (1111) \\ \hline
(4)      &  t &   \cdot         &  \cdot     & \cdot  \\
(13)    &  1      & \cdot     &  \cdot     & \cdot  \\
(22)    &  \cdot   & t          &  \cdot     &  \cdot  \\
(112)  &  \cdot      & 1           &    \cdot    &   \cdot \\
(31)    &   \cdot  &    \cdot      &   t      &   \cdot   \\
(121)  & \cdot       &  \cdot      &   1   & \cdot    \\
(211)  & \cdot       &  \cdot        & \cdot   & t  \\
(1111) & \cdot      &    \cdot      &  \cdot     &  1\\ \hline
 \multicolumn{5}{c}{} \\
\multicolumn{5}{c}{(13)\text{-Schur functions}} \\
\end{array} 
\quad \quad 
 \begin{array}{|c|cc|}
\hline 
          & (112) & (1111)\\ \hline
(4)      &  t^3   & \cdot      \\
(13)    &  t^2   & \cdot      \\
(22)    &  t       & \cdot      \\
(112)  &  1       & \cdot    \\
(31)    &  \cdot &   t^3   \\
(121)  & \cdot  &   t^2       \\
(211)  & \cdot   &   t       \\
(1111) & \cdot   &   1  \\ \hline
 \multicolumn{3}{c}{} \\
\multicolumn{3}{c}{(112)\text{-Schur functions}} \\
\end{array} \quad \quad 
\begin{array}{|c|cc|}
\hline 
          & (121)   & (1111) \\ \hline
(4)      &  t^4     & \cdot      \\
(13)    &  t^3     & \cdot      \\
(22)    & \cdot   & t^4      \\
(112)  &  \cdot   & t^3    \\
(31)    &  t         &   \cdot   \\
(121)  &  1        &   \cdot      \\
(211)  & \cdot   &   t       \\
(1111) & \cdot   &   1 \\ \hline 
\multicolumn{3}{c}{} \\
\multicolumn{3}{c}{(121)\text{-Schur functions}} \\
\end{array}
$$
\\
$$
\begin{array}{|c|cc|}
\hline 
          & (211)   & (1111) \\ \hline
(4)      &  t^5     & \cdot      \\
(13)    &  \cdot     & t^5      \\
(22)    & t^3   & \cdot      \\
(112)  &  \cdot   & t^3    \\
(31)    &  t^2         &   \cdot   \\
(121)  &  \cdot       &   t^2      \\
(211)  &  1   &   \cdot       \\
(1111) & \cdot   &   1 \\ \hline 
\multicolumn{3}{c}{} \\
\multicolumn{3}{c}{(211)\text{-Schur functions}} \\
\end{array}
$$
\subsection{Table of Macdonald polynomials in the $\gamma$-Schur basis in weight 4}  The tables below show examples of Theorem \ref{MacdonaldPositivity}.  The columns indicate
the index $\alpha$ of the function $\mathbf{H}_\alpha(A;q,t)$ .
$$
\begin{array}{|c|cccc|}
\hline
           & (31) & (121) & (211) & (1111) \\ \hline
(31)    &  1     &  t     & t^2       & t^3 \\
(121)  &  q^3  & 1    & q^3t^2  & t^2 \\
(211)  &  q^2  & q^2 &    1       &   t   \\ 
(1111) & q^5   & q^2 &   q^3   &   1 \\ \hline
\multicolumn{3}{c}{} \\
\multicolumn{5}{c}{\mathbf{H}\vert_{(31)}  \text{ in the } (31)\text{-Schur basis}}
 \end{array}
\quad \quad
\begin{array}{|c|cccc|}
\hline
           & (22) & (112) & (211) & (1111) \\ \hline
(22)    &  1     &  t     & t^3      & t^4 \\
(112)  &  q^3  & 1    & q^3t^3  & t^3 \\
(211)  &  q      & qt   &    1       &   t   \\ 
(1111) & q^4   & q   &   q^3   &   1 \\ \hline
\multicolumn{3}{c}{} \\
\multicolumn{5}{c}{\mathbf{H}\vert_{(22)}  \text{ in the } (22)\text{-Schur basis}}
 \end{array}
 $$
 \\
 $$
 \begin{array}{|c|cccc|}
\hline
           & (13) & (112) & (121)     & (1111) \\ \hline
(13)    &  1     &  t^2     & t^3       & t^5 \\
(112)  &  q^2  & 1        & q^2t^3  & t^3 \\
(121)  &  q      & qt^2   &    1       & t^2   \\ 
(1111) & q^3   & q   &   q^2         &   1 \\ \hline
\multicolumn{3}{c}{} \\
\multicolumn{5}{c}{\mathbf{H}\vert_{(13)}  \text{ in the } (13)\text{-Schur basis}}
 \end{array}
 $$

\end{document}